\batchmode
\documentclass [11pt,leqno]{article}
\usepackage{amssymb}
\newtheorem{theorem}{Theorem}

\newtheorem{proposition}[theorem]{Proposition}
\newtheorem{corollary}[theorem]{Corollary}
\newtheorem{definition}[theorem]{Definition}
\newtheorem{remark}[theorem]{Remark}

\newcommand{\X}{\mbox{${\partial}\kern-3pt {\cal C}$}}

\newcommand{\E}{\mbox{\,{{\rm E}\kern-7.6pt {\rule{.2mm}{2.5mm}}}\enskip\,}}
\newcommand{\pte}{\mbox{\,${{\sqsubset}\kern-1.8pt
           {\rule{.1mm}{2mm}}}\kern-6.8pt {\times}$\,}}

\title{ THE GEOMETRY OF A  BI-LAGRANGIAN MANIFOLD}
\author{Fernando Etayo\footnote{Departamento de Matem\'{a}ticas, Estad\'{\i}stica y
Computaci\'{o}n.
 Facultad de Ciencias.  Universidad de Cantabria.
 Avda. de los Castros, s/n, 39071 Santander, SPAIN.
 e-mail: etayof@unican.es},  Rafael Santamar\'{\i}a
 \footnote{Departamento de Matem\'{a}ticas. Escuela de Ingenier\'{\i}a Industrial e
 Inform\'{a}tica.
 Universidad de Le\'{o}n.
 Campus de Vegazana, 24071 Le\'{o}n, SPAIN. e-mail: demrss@unileon.es}\,  and Uju\'{e} R. Tr\'{\i}as \footnote{Departamento de
 Matem\'{a}ticas, Estad\'{\i}stica y
Computaci\'{o}n.
 Facultad de Ciencias.  Universidad de Cantabria.
 Avda. de los Castros, s/n, 39071 Santander, SPAIN.}}
\date{}
\begin{document}
\maketitle

\begin{abstract}
This is a survey on bi-Lagrangian manifolds, which are symplectic manifolds endowed with two transversal Lagrangian foliations. We also study the non-integrable case (i.e., a symplectic manifold endowed with two transversal Lagrangian distributions). We show that many different geometric structures  can be attached to these manifolds and we carefully analyse the  associated connections. Moreover, we introduce the problem of the intersection of two leaves, one of each foliation, through a point and show a lot of significative examples.
\end{abstract}


2000 Math. Sub. Class.: 53D05, 53D15, 53D17, 53C15, 53C05, 53C10, 53C12, 53C50, 57R30.

\section{Introduction}

A bi-Lagrangian manifold is a manifold $M$ endowed with three structures: a symplectic form $\omega$ , an integrable almost product structure $F$ (which defines two transversal equidimensional Lagrangian foliations) and a neutral metric $g$. In fact, two of the above structures determine the third one by means of the relation $\omega (X,Y)=g(FX,Y)$, for all $X,Y$ vector fields on $M$. So, bi-Lagrangian manifolds are in the intersection of three geometries: the symplectic, the almost product and the semi-Riemannian Geometry. Moreover, in a bi-Lagrangian manifold the Levi-Civita connection parallelizes the almost product and the symplectic structures. Bi-Lagrangian manifolds also named para-K\"{a}hler manifolds.

 In some sense, this work is a continuation of the survey \cite{CFG}, but there are important differences between both papers:

(1) We focus the attention  on symplectic aspects of the theory, i.e., on the geometry defined by $\omega$ and $F$ instead of that defined by $g$ and $F$.

(2)  We show many geometric structures that can be relationed with a bi-Lagrangian structure. This can allow to obtain more information about these manifolds. For example, we  study the the relation with special symplectic manifolds, Poisson structures and Lie algebroids.

(3) We present a complete study of the connections attached to bi-Lagrangian manifolds.

(4)We include some new problems, such as the study of the intersection of both leaves, one of each foliation, through a point.

Symplectic geometry is an active   topic of research, linking Differential and Algebraic Geometries, Algebraic Topology, Mathematical Physics and other fields. The reader can find several recent monographies about it, such as \cite{Ca} and \cite{MS}. Lagrangian foliations on symplectic manifolds are used in geometric
quantization. As is well known, the existence of a connection canonically
attached to a symplectic manifold is an important tool to obtain a
deformation quantization \cite{H}, \cite{Li}, \cite{Fe}. A
{\em bi-Lagrangian manifold} (i.e., a symplectic manifold endowed with two
transversal Lagrangian foliations) admits a canonical symplectic connection, which has been introduced by Hess in \cite{H}, and used by several
authors (e.g., \cite{Va} and \cite{Bo}). We shall name it the {\em bi-Lagrangian connection}. On the other hand, an {\em almost bi-Lagrangian manifold} (i.e., a symplectic manifold endowed with two
transversal Lagrangian distributions) also admits a canonical  connection, which is non-symmetric in general. We shall name it the {\em almost bi-Lagrangian connection}. The bi-Lagrangian and the almost bi-Lagrangian connection of a
bi-Lagrangian manifold coincide.

The present paper is a survey about the geometry of bi-Lagrangian manifolds and almost bi-Lagrangian manifolds. We choose the notion of bi-Lagrang\-ian structure in the general framework of Symplectic Geometry as the start point of this paper. Other geometric structures, such as those of para-Complex Geometry, will be introduced when necessary.  There are no complete proofs in the paper, but there are some ``Sketch of  Proof" and some elementary proofs. In particular, we prove the results linking different structures on a manifold. Examples are carefully explained.

\bigskip

An important remark is the following: there exists a different concept with the same name. Many authors name a 
 {\em bi-Lagrangian} manifold (resp. distribution) as a manifold (distribution) which is Lagrangian respect to two different symplectic structures (see, e.g., \cite{Sc}). This notion is relationed with that of a bi-symplectic and a bi-Hamiltonian structure, which depends on two different symplectic structures defined on the same manifold. We do not work with this definition in the present paper.

\bigskip

The paper is organized as follows:

In section 2 we present the  geometric properties of a
bi-Lagrangian manifold:  basic definitions, the bi-Lagrangian
connection,   the identity ``bi-La\-grang\-ian=para-K\"{a}hler", the associated $G$-structure, the different sectional curvatures, its automorphisms and symmetric bi-Lagrangian manifolds, the Poisson structure, the Lie algebroid associated to a 
bi-Lagrang\-ian structure, and the 3-web structures that one can
attach to such a manifold.

In section 3 we study connections on a bi-Lagrangian manifold, proving that  the well adapted, the Libermann, the natural, the  bi-Lagrangian, the almost bi-Lagrangian  and the Levi-Civita connections coincide.

In section 4 we study the holonomy of the leaves in the following sense: Let us consider a bi-Lagrangian manifold $M$, a point $p\in M$ and the leaves $L_{1}$ and $L_{2}$, one of each foliation, through the point $p$. In this section we shall obtain information about the number $N(p)$  of points in the intersection $L_{1}\cap L_{2}$. We shall say that $p$ has the trivial intersection property if $N(p)=1$. We shall prove that the following three concepts are independent: (1) trivial intersection property, (2)compactness of the manifold and (3) flatness of the canonical semi-Riemannian metric attached to $M$. We shall distinguish between the cases  ${\rm dim}(M)=2$ and  ${\rm dim}(M)>2$, because in the first one the manifold is also Lorentz.

We shall end the paper with some open questions.

\vspace{5mm}

All the manifolds through the paper will be assumed smooth. The
Lie algebra of vector fields of a manifold $M$ will be denoted as
${\cal X} (M)$. A Riemannian metric will be denoted as $G$,
whereas a semi-Riemannian metric of signature $(n,n)$ will be
denoted as $g$. On the other hand, automorphisms of ${\cal X} (M)$
of square $-id$ (resp. $id$) will be denoted as $J$ (resp. $F$ or
$P$).

\section{Symplectic and bi-Lagrangian manifolds}

In this section we shall present the basic definitions about Lagrangian structures and the connections attached to them and we shall obtain the first results.

\subsection{Lagrangian structures}

Let $(M,\omega )$ be a symplectic manifold, with $dim\: M=2n$. In the last
years, the following definitions have been introduced:
 A {\em Lagrangian distribution} on $M$ is a
$n$-dimensional distribution $\cal{D}$ such that $\omega (X,Y)=0$
for all vector fields $X,Y\in \cal{D}$. Such a Lagrangian
distribution is also called an {\em almost cotangent structure}
\cite{TS}.  A foliation $\cal{F}$ on $M$ is said a {\em Lagrangian
foliation} if its leaves are lagrangian submanifolds, i.e., each
leaf $N$ has $dim \: N=n$ and $\omega (X,Y)=0$, for every $X,Y$
tangent to $N$. A Lagrangian foliation is also called a {\em
polarization} \cite{AM} and an {\em integrable almost cotangent
structure}  \cite{TS}.

A symplectic manifold $(M,\omega )$ is said endowed with an {\em
almost bi-Lagrangian structure} (resp. {\em bi-Lagrangian
manifold}) if $M$ has two transversal Lagrangian distributions
(resp. involutive transversal Lagrangian distributions) ${\cal D}
_{1}$ and ${\cal D} _{2}$.  In this last case, the manifold is
endowed with two transversal foliations  ${\cal F} _{1},{\cal F}
_{2}$ whose tangent distributions ${\cal D} _{i}=T({\cal F}
_{i})$, $i=1,2$, define an almost bi-Lagrangian structure. We also
say that a bi-Lagrangian structure is an {\em integrable} almost
bi-Lagrangian structure.

An almost bi-Lagrangian manifold $(M,\omega ,{\cal D}
_{1},{\cal D} _{2})$ is an almost product manifold and then one can define a $(1,1)$ tensor field $F$
by $F\mid _{{\cal D} _{1}}=id$ and $F\mid
_{{\cal D} _{2}}=-id$, ${\cal D} _{i}$. Obviously , $F^{2}=I$ and the Nijenhuis $(1,2)$ tensor field $N_{F}$ vanishes iff both distributions are involutive. The projection over ${\cal D} _{1}$ (resp. ${\cal D} _{2}$) will be denoted by $\pi _{1}$ (resp. $\pi _{2}$): $\pi _{1}=\frac{I+F}{2}, \: \pi _{2}=\frac{I-F}{2}$. As the distribution ${\cal D} _{1}$ (resp. ${\cal D} _{2}$) is the eigenspace associated to the eigenvalue +1 (resp. -1) we also denote ${\cal D} _{1}=F^{+}$ and ${\cal D} _{2}=F^{-}$.

\newpage

\subsection{Bi-Lagrangian connections}

As is well known, a symplectic manifold $(M,\omega )$ admits
several symplectic connections (a {\em symplectic connection}
$\nabla $ is a torsionless connection parallelizing $\omega$), but
one needs additional assumptions to obtain a {\em canonical
connection} (see e.g. \cite{GRS}, where some sufficient conditions
are quoted. A symplectic manifold with a fixed symplectic
connection is called a {\em Fedosov manifold}). Bi-Lagrangian
manifolds admit a canonical connection, introduced by Hess in 1980
\cite{H}
 in a quite difficult way, that one can reduce to the
following expression (see also \cite{Va}, \cite{Bo}):  The {\em
bi-Lagrangian connection of a bi-Lagrangian manifold} is the
uni\-que symplectic connection $\nabla$ which parallelizes both
foliations ${\cal F} _{1}$ and ${\cal F} _{2}$, i.e., such that
$\nabla _{X}Y\in T({\cal F} _{i})$, for all vector field $X$ in
$M$ and all vector field $Y\in T({\cal F} _{i})$. If we define the
(1,1) tensor field $F$ by $F\mid _{{\cal D} _{1}}=id$ and $F\mid
_{{\cal D} _{2}}=-id$, ${\cal D} _{i}$ be the tangent distribution
to the foliation ${\cal F} _{i}$, it is easily shown that $\nabla
$ parallelizes both foliations ${\cal F} _{1}$ and ${\cal F} _{2}$
iff $\nabla F=0$. Then, a {\em bi-Lagrangian connection of a
bi-Lagrangian manifold} is the unique symmetric connection
satisfying $\nabla \omega=0$, $\nabla F=0$.

\vspace{3mm}

Observe that one cannot extend this definition to the case of an almost bi-Lagrangian manifold: $\nabla$ is torsionless because the Lagrangian distributions are involutive \cite[p. 569]{Va}. Nevertheless, one can give the following generalization of the above definition to the case of almost bi-Lagrangian manifolds (see \cite[p. 158]{H}:
 The {\em almost bi-Lagrangian connection} of an almost bi-Lagrangian manifold $(M,\omega ,{\cal D} _{1},{\cal D} _{2})$ is the unique connection $\nabla $ on $M$ which paralelizes $\omega $, ${\cal D} _{1}$ and ${\cal D} _{2}$ and verifies ${\rm Tor}_{\nabla }(X_{1},X_{2})=0,\, \forall X_{i}\in {\cal D} _{i}$, where ${\rm Tor}_{\nabla }$ denotes the torsion tensor of $\nabla $. As in the above case, we can say that the {\em almost bi-Lagrangian connection} of an almost bi-Lagrangian manifold $(M,\omega ,{\cal D} _{1},{\cal D} _{2})$ is the unique connection $\nabla $ on $M$ satisfying $\nabla \omega=0$, $\nabla F=0$ and ${\rm Tor}_{\nabla }(X_{1},X_{2})=0,\, \forall X_{i}\in {\cal D} _{i}$.

\vspace{3mm}

Obviously, the almost bi-Lagrangian connection of a bi-Lagrangian manifold is the bi-Lagrangian connection.

\subsection{Lagrangian distributions}

The existence of a Lagrangian distribution on a symplectic manifold implies the existence of  infinite Lagrangian distributions:

\begin{theorem} {\em \cite{ES2001}} Let $(M,\omega )$ be a symplectic manifold and let
$\cal{D}$ be a Lagrangian distribution. Then, $M$ admits infinitely many
different Lagrangian distributions.
\label{many}
\end{theorem}

{\em Sketch of  Proof.}

(1) Consider any Riemannian metric $G$ on $M$ and let $\cal{D}^{\perp }$ the $G$-orthogonal distribution. One can prove that $\cal{D}^{\perp }$  is a Lagrangian distribution.

(2) Define the almost Hermitian structure $(M,J,G)$ associated to $(M,\omega )$, i.e., $\omega (X,Y)=G(JX,Y)$, for all vector fields $X,Y$ on $M$.

(3) Define the endomorphism of vector fields $F$ given by $F(X)=X$, if $X\in \cal{D}$, and $F(X)=-X$, if $X\in \cal{D}^{\perp }$. Then, one proves that $G(FX,FY)=G(X,Y)$, for all vector fields $X,Y$. Thus, $(M,F,G)$ is a Riemannian almost product manifold.

(4) Let $\alpha , \beta \in \mathbb{R}$ such that $\alpha ^{2}+\beta ^{2}=1$. Then $F_{(\alpha ,\beta )}$, given by $F_{(\alpha ,\beta )}(X)=\alpha F(X)+\beta JF(X)$,  is an almost product structure whose eigenspaces define Lagrangian distributions.$\Box$

\begin{remark}
\label{plucker}

The above result shows that  a symplectic manifold endowed with a Lagrangian distribution admits infinitely many distributions. In some cases, one can determine them. For example, if $M$ is the real plane then every straight line is a Lagrangian submanifold. Now, let us consider the space $\mathbb{R}^{4}$ endowed with the symplectic structure given by
$\omega =\left( \begin{array}{cc} 0 &I_{2} \\ -I_{2} & 0\end{array} \right) $, where $I_{2}$ denotes the $2\times 2$ identity matrix. We want to determine its  Lagrangian planes. This will be made by means of Pl\"{u}cker coordinates in the Grassmann manifold  $G(2,4)$  of vector planes in ${\mathbb R}^{4}$. We use the notation of {\rm \cite{E2003}}.

One can immerse $G(2,4)$ in $P_{5}({\mathbb R})$ by using Pl\"{u}cker coordinates: if the vectors $a=(a_{1},a_{2},a_{3},a_{4})$ and  $b=(b_{1},b_{2},b_{3},b_{4})$ define a basis of the plane $\Pi $ , then one can define the numbers $p_{ij}=\left| \begin{array}{cc} a_{i} & a_{j} \\ b_{i} & b_{j} \end{array} \right| $, with $1\leq i<j\leq4$ and one can associate to the plane $\Pi $ the homogeneous coordinates  $(p_{12}:p_{13}:p_{14}:p_{23}:p_{24}:p_{34})\in P_{5}({\mathbb R})$.  As is well known, under a change of the basis in the plane $\Pi $, the Pl\"{u}cker coordinates are multiplied by the determinant of the transformation matrix and
 $G(2,4)$ can be identified with the quadric  $\{ p_{12}p_{34}-p_{13}p_{24}+p_{14}p_{23}=0\} \subset P_{5}({\mathbb R})$.

 Let $\Pi$ be a plane generated by two independent vectors $v=(a,b,c,d)$ and $w=(\alpha ,\beta , \gamma ,\delta )$. Then , $\omega _{\Pi }=0$ iff  $\omega (v,w)=0$, i.e., $p_{13}+p_{24}=0$, which shows that the set of Lagrangian planes is a
3-dimensional manifold (because ${\rm dim}\, G(2,4)=4$). Moreover, as is well  known, Lagrangian planes are totally real planes when one considers the standard K\"{a}hler structure ${\mathbb C}^{2}={\mathbb R}^{4}$ (see {\rm \cite[Proposition 12]{E2003}}).

\end{remark}

On the other hand, one can easily prove that the bi-Lagrangian connection associated to certain bi-Lagrangian structures defined in a K\"{a}hler manifold is the Levi-Civita connection of the Riemannian metric of the manifold:

\begin{theorem}{\em \cite{ES2001}} Let ${\cal F}$ be a Lagrangian foliation in a K\"{a}hler
manifold $(M,J,G)$, such that the Levi-Civita connection $\nabla$ of $G$
parallelizes the foliation. Then:

(1) The orthogonal distribution ${\cal D}^{\perp }=(T{\cal F})^{\perp }$ is
parallel with respect to $\nabla$.

(2) ${\cal D}^{\perp }$ is involutive, and then $M$ is a bi-Lagrangian
manifold.

(3) $\nabla$ is the bi-Lagrangian  connection associated to the bi-Lagrangian
structure.

(4) All the distributions obtained in  Theorem \ref{many} are involutive.

\end{theorem}

\begin{remark}
If $(M,J,G)$ is a  K\"{a}hler manifold of dimension $2n$ then Lagrangian submanifolds coincide with totally real 
$n$-dimensional submanifolds, i.e., $N$ is a Lagrangian submanifold iff $J(T_{p}N)=T^{\perp }_{p}M$, for all $p\in N$. In the survey \cite{Chen2001} one can find a complete study of Lagrangian submanifolds of  K\"{a}hler manifolds.
\end{remark}

\begin{remark}
\label{cuadrica}

 The following example shows that a K\"{a}hler manifold admitting two transversal foliations may be a no bi-Lagrangian manifold. Let us consider the {\em complex quadric} $Q_{2}(\mathbb{C})\subset P_{3}(\mathbb{C})$ defined $\{ z_{0}^{2}+z_{1}^{2}+z_{2}^{2}+z_{3}^{2}=0\} $, where $z_{i}$ denotes the homogeneous coordinates in $P_{3}(\mathbb{C})$. The geometric properties of this manifold are well known (see, e.g. \cite[[pp.478-480]{GH}, \cite[p. 13]{Ha}, \cite[Example 10.6, chap. 11]{KN}  ). In particular, $Q_{2}(\mathbb{C})$ is a K\"{a}hler Einstein manifold which is isomorphic to $P_{1}(\mathbb{C})\times P_{1}(\mathbb{C})$   as K\"{a}hler  manifolds. The quadric contains two families of complex projective lines in such a way that two different lines of the same family do not intersect and a line of a family meets every line of the other family in exactly one point.

Then, the quadric $Q_{2}(\mathbb{C})$ has two transversal foliations, but they are not Lagrangian because the complex structure on $Q_{2}(\mathbb{C})\thickapprox P_{1}(\mathbb{C})\times P_{1}(\mathbb{C})$ restricts to each leave of both foliations. Then if $(Q_{2}(\mathbb{C}),J,G)$ denotes the K\"{a}hler structure of the quadric,  one has $\omega (X,JX)=G(JX,JX)=G(X,X)\neq 0$, for any $X\neq 0$ vector tangent to a leaf. Then, the leaves are not Lagrangian submanifolds and $Q_{2}(\mathbb{C})$ is not a bi-Lagrangian manifold.

\end{remark}

\subsection{Para-K\"{a}hler structures}

The main result in this section shows that ``bi-Lagrangian=para-K\"{a}hler". First, we show that the bi-Lagrangian connection is always the Levi-Civita connection of a neutral metric of signature $(n,n)$.

\begin{theorem} {\em  \cite{ES2001, ES2001b}}  Let $(M,\omega ,{\cal F} _{1},{\cal F} _{2})$ be a
bi-Lagrangian manifold. Then $M$ admits a canonical neutral metric $g$ whose
Levi-Civita connection coincides with the bi-Lagrangian connection of the
bi-Lagrangian manifold.
\label{metrica}
\end{theorem}

{\em Sketch of  Proof.}

(1) Let $(M,\omega ,{\cal F} _{1},{\cal F} _{2})$ be a
bi-Lagrangian manifold. We define the (1,1) tensor field $F$ by $F\mid
_{{\cal D} _{1}}=id$ and $F\mid _{{\cal D} _{2}}=-id$, ${\cal D} _{i}$ be
the tangent distribution to the foliation ${\cal F} _{i}$, and the map $g$
which
applies two vector fields $X,Y\in {\cal X}(M)$ to $g(X,Y)=\omega (FX,Y)$.
Then $g$ is a neutral metric.

(2) Let $\nabla$ be the bi-Lagrangian connection of $(M,\omega ,{\cal F}
_{1},{\cal
F} _{2})$; then $\nabla$ is a torsionless connection. One can prove that  $\nabla g=0$, by using that $\nabla $
 parallelizes  $\omega$ and both foliations (or equivalently, $\nabla F=0$)$\Box$

\bigskip

Let us remember some basic definitions of Para-Complex Geometry (see the foundational works of Rashevskij \cite{R}  and Libermann \cite{Lb} and the survey \cite{CFG} of Cruceanu, Fortuny and  Gadea and the more than 100 references therein).  An {\em almost para-K\"{a}hler manifold} $(M,F,g)$ is a manifold endowed with a (1,1) tensor field $F$ satisfying $F^2=$id and a neutral metric $g$  such that $g(FX,FY)=-g(X,Y)$, for all vector fields  $X,Y$ in $M$, in such a way that the symplectic form $\omega $ defined by $\omega (X,Y)=g(FX,Y)$ is closed. A {\em para-K\"{a}hler manifold} is an almost para--K\"{a}hler manifold
$(M,F,g)$ such that $\nabla F=0$, $\nabla$ being the Levi-Civita connection
of $g$. Equivalently, both distributions $F^{+}$ and $F^{-}$, associated to the eigenvalues +1 and -1 of $F$  are involutive
and $\omega $ is closed.

Then, as a consequence of Theorem \ref{metrica}, one can easily prove the following

\begin{proposition}  {\em \cite{ES2001, ES2001b}} Let $M$ be a manifold.

a) There exists a bijection between
almost bi-Lagrangian structures on $M$ and almost para-K\"{a}hler structures on $M$.

b) There exists bijection between  bi-Lagrangian structures on $M$ and para-K\"{a}hler structures on $M$.

\label{equivalencia}
\end{proposition}

The following table shows the different names of these manifolds,
when one uses (s) "symplectic" or (p) "para-complex" terminology.
Remember that one has three objects: a neutral metric $g$, an
almost symplectic form $\omega $ and an almost product structure$F$
with the relations: $g(FX,FY)=-g(X,Y); \: \omega (X,Y)
\\=g(FX,Y);\: \omega (FX,FY)=-\omega (X,Y)$, for all vector fileds
$X,Y$.

\begin{center}
\begin{tabular}{|c|c|c|}
\hline
two distributions & one distribution  & two foliations\\
  & and one foliation &  \\
\hline
(p) almost para-Hermitian & (p) 1-para-Hermitian & (p) para-Hermitian \\
\hline
  & (p) 1-para--K\"{a}hler &\\
  & $i_{X}i_{Y}d\omega =0$ & \\
  & for $X,Y$ in the foliation  & \\
\hline
(p) almost para-K\"{a}hler & & (p) para--K\"{a}hler\\
(s) almost bi-Lagrangian & & (s) bi-Lagrangian\\
$ d\omega =0$ &   &$ d\omega =0$ \\
\hline
\end{tabular}
\end{center}

The geometry of 1-para-Hermitian  and 1-para--K\"{a}hler manifolds have been studied in \cite{CE}. As is well known, the cotangent bundle $T^{*}M$ of a manifold $M$ is endowed with a canonical symplectic structure and the vertical distribution is a Lagrangian distribution, i.e., the cotangent bundle has an almost cotangent structure \cite{TS}. If the manifold is endowed with a linear connection then its cotangent bundle is an 1-para--K\"{a}hler manifold.

Some results about almost symplectic manifolds endowed with two Lagrangian distributions or foliations (i.e., about almost para-Hermitian and para-Hermitian manifolds) have been studied in the paper \cite{GRS}. Nevertheless, these manifolds have no ``symplectic denomination".

\subsection{G-structure and topological obstructions}

Taking into account the results of the above section one has:

\begin{proposition} {\em \cite{ES1999}} The structure group of almost bi-Lagrangian manifolds is the paraunitary group 

\[U(n , {\bf A}) = \left\{ \left(
\begin{array}{ll}
A & 0 \\
0 & ^t A^{-1}
\end{array}
\right)
\colon \; 
A \in GL (n; \mathbb R) \right\}.\]

\end{proposition}

Moreover, one has the following information about the paraunitary Lie algebra

\[ {\mathfrak u}(n , {\bf A})= \left\{ \left(
\begin{array}{ll}
A & 0 \\
0 & -A^{t}
\end{array}
\right)
\colon \; 
A \in \mathfrak{gl} (n; \mathbb R)\right\} .\]

\noindent which will be useful to us in the study of connections attached to such a manifold.

\begin{proposition}{\em \cite{ES1999}}  The Lie algebra ${\mathfrak u}(n,{\bf A})$ is invariant under matrix transposition and its first prolongation vanishes.
\label{teor:buenascondiciones}
\end{proposition}

First topological obstructions can be found taking into account that bi-Lagrangian manifolds are symplectic (then they are orientable and even-dimen\-sional) and para-Hermitian (and then, the tangent bundle can be decomposed as the Whitney sum of two subbundles with the same rank). In \cite{G} one can find explicit obstructions by means of the Euler characteristic class. In particular, the following manifolds cannot admit an almost bi-Lagrangian structure: the spheres $S^{n}$,  the real projective spaces $P_{n}(\mathbb{R})$, the complex projective spaces $P_{n}(\mathbb{C})$, the quaternionic projective spaces $P_{n}(\mathbb{H})$ and the product of spheres $S^{n}\times S^{m}$, with $n\neq m$.

In Section 4 we shall give a collection of examples in both situations, compact and no compact.

\subsection{Metric, para-holomorphic and symplectic sectional curvature}

As we have seen, a bi-Lagrangian manifold has a canonical neutral metric $g$. Then, one can study the sectional curvature defined by this metric associated to non-$g$-isotropic planes. When one restricts to $F$-invariant planes, then one say that it is the {\em para-holomorphic sectional curvature}. When it is constant over all the $F$-invariant
non-$g$-isotropic planes one says that the manifold is a {\em para-holomorphic space form}. A classification of non-flat
para-K\"{a}hler space forms was obtained in \cite{GM1995}. Unlike the complex case, constant para-holomorphic sectional curvature $c\neq 0$ does not imply bounded sectional curvature (see \cite{GM1}, \cite{E3}): sectional curvature over non-$F$-invariant planes runs over all the real line $\mathbb{R}$.

On the other hand, para-K\"{a}hler space forms are {\em Osserman manifolds}, i.e., the eigenvalues of the Jacobi operator are constant (see \cite{GVV}; in that paper one can find examples of neutral manifolds which are nonsymmetric. A complete monograph about this topic is \cite{GKV}).

In the paper \cite{GRS} the authors introduce the notion of sectional curvature of a Fedosov manifold, i.e., a symplectic manifold endowed with a fixed symplectic connection. As bi-Lagrangian manifolds are Fedosov manifolds then the following question is natural: coincide the symplectic sectional curvature with the metric sectional curvature? 

Let us denote by $R^{\omega }$ (resp. $R^{g}$) the Fedosov curvature (resp. the semi-Riemannian curvature), which are defined as

\[
R^{\omega} (X,Y,Z,W)=\omega(X,R(Z,W)Y), \: R^{g}(X,R(Z,W)Y),
\forall X,Y,Z,W \in {\cal X}(M),
\]

\noindent  $R$ being the curvature tensor field  respect to the connection of the manifold. In the first case this connection is the bi-Lagrangian connection and, in the second one, it is the Levi-Civita connection, but they coincide (see Theorem \ref{metrica}) and then $R$ is the same $(1,3)$ curvature tensor field. 

\begin{proposition} With the above notation, if $M$ is a bi-Lagrangian manifold, then

\begin{eqnarray*}
R^{\omega } (X,Y,Z,W)&=& \omega (\pi _{2}X , R(Z,W) \pi _{2}Y)+ \omega (\pi _{2}X,
R(Z,W) \pi _{1}Y),\\
\noalign{\medskip} R^{g} (X,Y,Z,W)&=& \omega (\pi _{2}X , R(Z,W) \pi _{2}Y)-
\omega (\pi _{2}X, R(Z,W) \pi _{1}Y).
\end{eqnarray*}

\label{comparacion}
\end{proposition}

{\em Proof.} It follows from a direct computation. For $R^{\omega }$ we have:

\begin{eqnarray*}
R^{\omega} (X,Y,Z,W)  = \omega(X,R(Z,W)Y)\\ = \omega (\pi _{1}X +\pi _{2}X,
R(Z,W)(\pi _{1}Y + \pi _{2}Y))\\
= \omega (\pi _{1}X, R(Z,W)\pi _{1}Y)+ \omega(\pi _{1}X , R(Z,W) \pi _{2}Y)\\ +
\omega
(\pi _{2}X , R(Z,W) \pi _{1}Y) + \omega (\pi _{2}X , R(Z,W) \pi _{1}Y)\\
= \omega (\pi _{2}X , R(Z,W) \pi _{2}Y)+ \omega (\pi _{2}X,
R(Z,W) \pi _{1}Y),
\end{eqnarray*}

\noindent taking into account that $\nabla F=0$, thus preserving both distributions, and $\omega $ vanishes over them. A similar calculation proves the result for $R^{g}\: \Box$

\subsection{Automorphisms and symmetric bi-Lagrangian manifolds}

Let $(M,\omega ,{\cal F} _{1},{\cal F}_{2})$ be a bi-Lagrangian manifold. An important problem consists on the determination of the group of automorphisms of $M$ which preserve the bi-Lagrangian structure. We begin with the definition of an automorphism.

\begin{definition} {\em \cite{K2}} 
\label{automorfismos}
Let $(M,\omega ,{\cal F} _{1},{\cal F}_{2})$ be a bi-Lagrangian manifold.

(a) A {\em paracomplex automorphism} of $M$ is a diffeomorphism $\varphi $ of $M$ which leaves the leaves of both foliations ${\cal F} _{1},{\cal F}_{2}$  invariant, i.e., $\varphi _{*}\circ F=F\circ \varphi _{*}$.

(b) A {\em paracomplex isometry} of $M$ is a paracomplex automorphism leaving $\omega $ invariant, i.e., $\varphi ^{*}(\omega )=\omega $.

\end{definition}

\begin{proposition}
\label{isometriag}
Let $(M,\omega ,{\cal F} _{1},{\cal F}_{2})$ be a bi-Lagrangian manifold.  Let $\varphi :M\to M$ is a paracomplex automorphism. Then, $\varphi $ is a paracomplex isometry iff it is a $g$-isometry, $g$ being the canonical neutral metric attached to $M$ by Theorem \ref{metrica}.

\end{proposition}

{\em Proof.}  The result follows from a direct calculation:  let $X,Y$ be vector fields on $M$.

$\Rightarrow )\; g(\varphi _{*}X,\varphi _{*}Y)= g(FF\varphi _{*}X,\varphi _{*}Y)= \omega (F\varphi _{*}X,\varphi _{*}Y)= \omega (\varphi _{*}FX,\varphi _{*}Y)=  \omega ^{*}(FX,Y)= \omega (FX,Y)= g(FFX,Y)= g(X,Y).$

$\Leftarrow )\;   \omega ^{*}(X,Y)= \omega (\varphi _{*}X,\varphi _{*}Y)  = g(F\varphi _{*}X,\varphi _{*}Y) = 
g(\varphi _{*}FX,\varphi _{*}Y) =\\  g(FX,Y) = \omega (X,Y).   $

The result is proved$\Box$

Observe that one can conclude then that the group of paracomplex isometries coincides with the group of  
 $g$-isometries of the manifold, and then, it is a Lie group of dimension $\leq {\rm dim}(M)$ (cfr. \cite[pp. 255  and 258]{O'N}). The group of paracomplex automorphisms is not in general a  finite-dimensional Lie group \cite[p. 536]{K2}.

\begin{definition} {\em \cite{K2}} A symmetric space $M=G/H$ is called a {\em bi-Lagrangian symmetric space} if $M$ admits a $G$-invariant bi-Lagrangian structure, i.e., if $M$ admits a bi-Lagrangian structure $(M,\omega ,{\cal F} _{1},{\cal F}_{2})$ such that $G$ acts on $M$ as paracomplex isometries.

\end{definition}

Bi-Lagrangian symmetric spaces with $G$ a semisimple Lie groups have been classified by Kaneyuki and Kozai \cite{KK}. On the other hand, the problem of detemining the group paracomplex automorphisms of a bi-Lagrangian symmetric space remains open in general, although there are significative advances (see \cite{K2} and the references therein). And the general problem is to determine the group of paracomplex automorphisms of a bi-Lagrangian manifold.

\subsection{Special complex manifolds and Poisson structures}

In last years several papers about special complex, symplectic and K\"{a}hler manifolds have been published. We shall follow the notation of \cite{ACD}. A {\em special complex manifold} $(M,J,\nabla )$ is a complex manifold $(M,J)$ together with a flat torsionfree connection $\nabla$ such that $\nabla J=0$. A {\em special symplectic manifold} $(M,J,\nabla ,\omega )$ is a  special complex manifold $(M,J,\nabla )$ together with a $\nabla $-parallel symplectic structure $\omega$. These manifolds can be immersed in the cotangent bundle $T^{*}\mathbb{C}^{n}$, $n$ being the complex dimension of $M$, in such a way that $\omega =2\sum dx^{i}\wedge dy_{i}$, where $\{ x^{1},\ldots ,x^{n}, y_{1},\ldots ,y_{n}\}$ is real coordinate system near each point of $M$.

One can try to obtain a similar theory for bi-Lagrangian manifolds. In this case, we have the symplectic form $\omega$ and the almost product $F$ (instead of the complex structure $J$). The following result of  Boyom  \cite[Theorem 1.2.1]{Bo}  can be  viewed as the first result in this topic :

\begin{theorem}
\label{Boyom}

Let $(M,\omega ,{\cal F} _{1},{\cal F}_{2})$ be a bi-Lagrangian manifold and let $\nabla$ be the bi-Lagrangian connection. Then the following conditions are equivalent:

(1) $\nabla$ is a flat connection.

(2) For every point in $M$ there exists a coordinate system $\{ x^{1},\ldots ,x^{n}, y_{1},\ldots ,y_{n}\}$ such that: (a) 
 $\{ x^{i},x^{j}\} =\{ y_{i},y_{j}\} =0$ and $\{ x^{i},y_{j}\} =\delta ^{i}_{j}$; (b) the distribution ${\cal D} _{1}=T({\cal F} _{1})$ (resp. ${\cal D} _{2}=T({\cal F} _{2})$) is locally generated by the Hamiltonian vector fields $X_{x^{1}},\ldots ,X_{x^{n}}$ (resp. $X_{y_{1}},\ldots ,X_{y_{n}}$), where $\{ \, ,\, \}$ denotes the Poisson bracket defined by the symplectic structure.

\end{theorem}

\subsection{The Lie algebroid associated to a bi-Lagrangian structure}

In last years the notion of {\em Lie algebroid} has been developed,   providing a general framework for different notions such as Lie algebras, bundles of Lie algebras, tangent bundles, etc., and relating several topics such as Poisson geometry, theory of connections, structures on manifolds, etc., (see \cite{CW, Fer} for  global expositions). A {\em Lie algebroid} over $M$ is a vector bundle $\pi :E\to M$, endowed with a Lie algebra structure $\{ -,-\}$ on $\Gamma (\pi )$, together a bundle map $\rho :E\to TM$ (called the {\em anchor}) such that

1) The map $\Gamma (\pi )\to \Gamma (TM)=\mathfrak{X}(M)$ is a Lie algebra homomorphism.

2) For any function $f$ on $M$ and any sections $s,s'\in \Gamma (\pi )$ the following identity holds:

$$\{ s,fs'\} =f\{ s,s'\} + (\rho \circ s)(f) s'$$

Following \cite{CL} we shall show that a bi-Lagrangian manifold can be endowed with a Lie algebroid. If $M$ is a manifold endowed with a (1,1) tensor field $F$ with vanishing Nijenhuis torsion $N_{F}=0$, then one can define a new bracket  by

$$\{ X,Y\} = [X,Y]_{F}=[F(X),Y]+[X,F(Y)]-F([X,Y])                      $$

\noindent $X,Y$ being vector fields over $M$. Moreover, one has $F([X,Y]_{F})=[F(X),F(Y)]$, which proves the first property  above, and $[X,fY]_{F}=f[X,Y]_{F} + F(X)(f)Y$, where $f$ is a function on $M$, which proves the second one, thus showing that the tangent bundle $TM$ becomes a Lie algebroid with bracket $[-,-]_{F}$ and anchor map $F:TM\to TM$.

Moreover, Fernandes \cite{Fer} defined the notion of Lie algebroid connection; in the present situation a
 pseudo-connection whose fundamental tensor field is $F$ is a Lie algebroid connection (see \cite{E1993} for a survey on pseudo-connections).

In our case of a bi-Lagrangian manifold,  taking  $F$ as the (1,1) tensor field corresponding to its
para-K\"{a}hler structure, which verifies $N_{F}=0$, then we can endow $M$ with  the Lie algebroid  with bracket $[-,-]_{F}$ and anchor map $F:TM\to TM$.

\subsection{Bi-Lagrangian manifolds and 3-webs}

We shall show that every bi-Lagrangian manifold can be endowed with a metric 3-web structure, by means of a  Riemannian metric. For the sake of completeness of this survey we first remember the basic features of 3-web structures.

Blaschke introduced the notion of 3-web on a surface as three families of curves on the surface such that any two of three curves through any point of the surface are always transversal. This notion has been extended to three distributions on an even dimensional manifold:  A {\em 3-web} on a manifold $M$ is given by  three equidimensional supplementary distributions: for all $x\in M$, the tangent space of $M$ at $x$ is decomposed as $T_{x}M=V_{1}(x)\oplus
V_{2}(x)=V_{1}(x)\oplus V_{3}(x)=V_{2}(x)\oplus V_{3}(x)$, $V_{1},V_{2},V_{3}$ being the distributions.

There exists an alternative presentation by means of polynomic structures on the manifolds: Cruceanu introduced in \cite{C} the notion of an {\em almost biparacomplex manifold} in the
following way: An {\em almost biparacomplex structure} on a manifold
$M$ is given by two anticommutative almost product structures $F$ and $P$,
i.e., two tensor fields $F$ and $P$ of type $(1,1)$  verifying
$F^{2}=P^{2}=1,\: F\circ P+P\circ F=0$. Then, there are four equidimensional
and supplementary distributions, defined by the eigenspaces associated to
$+1$ and $-1$ of the automorphisms $F$ and $P$ (namely $F^{+}$, $F^{-}$,
$P^{+}$, $P^{-}$). In particular, $M$ has even dimension, $F$ and $P$ are
almost paracomplex structures (because $dim F^{+}=dimF^{-},\:
dimP^{+}=dimP^{-}$) and $F$ (resp. $P$) is an isomorphism between $P^{+}$
and $P^{-}$ (resp. between $F^{+}$ and $F^{-}$).

Then, we can state

\begin{proposition} {\em $\cite{C}$}  A  manifold $M$ is endowed with a 3-web iff it is endowed with an almost biparacomplex structure.
\end{proposition}

{\em Sketch of  Proof.}
If $F$ is the almost product structure given by $F^{+}=V_{1}, F^{-}=V_{2}$ and $P$ is the almost
product structure given by $P^{+}=V_{3}, P^{-}=F(V_{3})=V_{4}$, one easily
check that $(M,F,P)$ is an almost biparacomplex manifold.$\Box$

\vspace{5mm}

Moreover, if $(M,F,P)$ is an almost  biparacomplex manifold, then
one can consider $J=P\circ F$, which is an {\em almost complex
structure} on $M$. So, an almost biparacomplex manifold is an
even-dimensional orientable manifold which has two  almost product
structures and one almost complex one.  One also can define an
{\em almost tangent} structure $K$ given by: $K(X)=P(X)$, if $X\in
F^{+}$ and $K(X)=0$, if $X\in F^{-}$.

\vspace{5mm}

One of us has introduced the following  metrics adapted to an almost  biparacomplex structure.

\begin{definition} {\em (see \cite{S})}.  Let $(M,F,P)$ be a biparacomplex
manifold and let $g$ be
a pseudo-Riemannian metric on $M$. Then, $(M,F,P,g)$ is said a $(\varepsilon
_{1},\varepsilon _{2})$ {\em pseudo-Riemannian almost biparacomplex manifold}, where
$\varepsilon _{1},\varepsilon _{2}\in \{+,-\}$
according to the following relations:

$$g(FX,FY)=\varepsilon _{1}g(X,Y);\:  g(PX,PY)=\varepsilon _{2}g(X,Y).$$
\end{definition}

Finally, an almost biparacomplex manifold $(M,F,P)$ is said a {\em biparacomplex manifold} if the distributions $F^+, F^-, P^+, P^-$ are involutive (or equivalently, if $N_{F}=N_{P}=0$, $N$ being the Nijenhuis tensor field).

\vspace{3mm}

By direct computations one can prove the following

\begin{proposition}
\label{web}
{\em \cite{ES2000}} Let $(M,\omega ,{\cal D} _{1},{\cal D} _{2})$ be an almost
bi-Lagrangian manifold and let $(M,F,g)$ be its associated almost para-K\"{a}hler structure.   For each Riemannian metric $G$ such that ${\cal D} _{1}$ and ${\cal D} _{2}$ are $G$-orthogonal, we define the almost complex structure $J$ associated to $G$ and $\omega $ (i.e., $\omega (X,Y)=G(JX,Y)$). Then:

(1) $(M,F,P=J\circ F)$ is an almost  biparacomplex manifold;

(2) $(M,J,g)$ is a Norden manifold;

(3) $(M,F,G)$ is a Riemannian almost product manifold;

(4) $(M,F,P,g)$ is a $(-,+)$ pseudo-Riemannian almost biparacomplex manifold;

(5) $(M,F,P,G)$ is a $(+,+)$ Riemannian almost  biparacomplex  manifold.

\end{proposition}

Such a metric always exists: if $H$ is any Riemannian metric on $M$, then $G$ given by $G(X,Y)=H(X,Y)+H(FX,FY)$ then  $(M,G,F)$ is an almost product manifold.

\section{Connections on an almost bi-Lagrangian manifold}

In the above sections we have introduced the almost bi-Lagrangian connection, which is not the Levi-Civita connection of the almost para-K\"{a}hler structure attached to the almost bi-Lagrangian one, and the bi-Lagrangian connection, which can be defined only in the integrable case (in this case, it  coincides with the Levi-Civita connection). We shall show that there exists other  connections  attached to an almost
bi-Lagrangian manifold. We shall prove that all of them coincide when the structure is integrable.  First, we shall
remember some results about functorial connections.

\subsection{Functorial connections}

We shall follow the notation of \cite{ES1999}. A {\em functorial connection} associated to a $G$-structure is, roughly speaking, a reducible connection which is natural with respect to the isomorphisms of the $G$-structure. Such connections are useful in the study of the integrability  of the $G$-structure and the calculus of the differential invariants of the
$G$-structure.  Moreover, the non-existence of a functorial connection makes the construction of differential invariants extremely difficult. It is the case of the symplectic and conformal geometries. Symplectic manifolds do not admit a functorial connection because the first prolongation of the Lie algebra associated to its structure group does not vanish. Nevertheless, as we shall show in this section, bi-Lagrangian manifolds admit functorial connections. Moreover, we shall show that one can define four functorial connections on an almost bi-Lagrangian manifold (the well adapted, the Libermann, the bi-Lagrangian  and the Levi-Civita connections), which coincide if the manifold is bi-Lagrangian.

\vspace{5mm}

We shall need basic results about functorial connections. Let $\pi\colon F(M) \to M$ be the bundle of linear frames and let  $\pi \colon P \to M$ be  a $G$-structure over  $M$.

\begin{definition}{\em (\cite[Def. 2.2]{Jaime})} A functorial connection is an assignment $\sigma \mapsto \nabla (\sigma)$, that associates a linear connection  $\nabla (\sigma)$ over $M$ to each section $\sigma$ of the bundle $F(M)/G$, satisfying the following properties:

i)  $\nabla (\sigma)$ is reducible to  the subbundle $P_{\sigma}$; i.e.,
$\nabla (\sigma)$ is adapted to $\sigma$.

 ii) For every diffeomorphism $f$ of $M$,  $\nabla (f\cdot \sigma)= f
      \cdot \nabla (\sigma)$, where  $f \cdot \nabla (\sigma)$ is the connection image of  $\nabla (\sigma)$ by $f$  in the subbundle   $P_{f \cdot \sigma} = {\widetilde f}(P_{\sigma})$.

iii) $\nabla (\sigma)$  depends continuously on $\sigma$ with respect  to the ${\cal C}^{\infty}$ topologies of the spaces of sections  of the classifying bundle and of the bundle of linear connections.

The last condition is equivalent to: there exists an integer $r \geq 0$ such that $\nabla (\sigma) (x)$ only depends on $j^r_x \sigma$, for every point $x \in M$.
\label{teor:deffuntorial}
\end{definition}

The vanishing of the first prolongation of the Lie algebra of $G$ provides an obstruction to the existence of functorial connections attached to the $G$-structure.

\begin{theorem}{\em (\cite[Teor. 2.1]{Jaime})}
If a functorial connection exists for  $G$-structu\-res, then the first prolongation of the Lie algebra of $G$ must vanish.
\label{teor:necesaria}
\end{theorem}

The two following results show the interest of this theory: one can obtain sufficient conditions for the existence of the {\em well adapted connection}, which is a functorial connection which mesure the integrability of the $G$-structure: it is integrable iff the tensors of torsion and curvature of the well adapted connection vanishes. In  \cite{Valdes}, we find

\begin{theorem} \begin{em} \cite[Teor. 1.1]{Valdes} \end{em}
The following two assertions are equivalent:
 \newcounter{Item30}
\begin{list}{ \roman{Item30})} {\usecounter{Item30} }
      \item For every  $G$-structure $ P \to M$, there exists a unique connection    $\nabla $
      adapted to the $G$-structure such that, for every endomorphism $S$ given by a section of the adjoint
bundle of $P$ and every vector field  $X \in {\cal X}(M)$,  one has
      \[ {\rm trace} (S \circ i_X \circ {\rm Tor}_{\nabla }) = 0 . \]
Moreover this connection only depends on the first contact  of the $G$-structure.

    \item If  $T \in {\rm Hom} (\mathbb{R}^{n} , \mathfrak g)$ satisfies that  $i_v \circ {\rm alt} (T)
    \in \mathfrak g^{\perp}$ for any   $v \in \mathbb{R}^{n}$,  then $T = 0 $,
    where $\mathfrak g$ is the Lie algebra of $G$,   $\mathfrak g^{\perp}$
    is the orthogonal subspace to $\mathfrak g$ in $\mathfrak{gl}(n,\mathbb{R})$
    with respect to the Killing-Cartan metric, and ${\rm alt} (T) (u,v) = T(u)v
    - T(v)u , \quad \forall u , v \in \mathbb{R}^{n}$.
\end{list}

The connection  $\nabla $ (if there exists) is called  the {\em well-adapted connection} associated to the  $G$-structure.
\end{theorem}

Taking into account Proposition \ref{teor:buenascondiciones} we obtain 

\begin{corollary}  One can define the well-adapted connection on every almost bi-Lagrangian manifold.
\label{existen}
\end{corollary}

\begin{theorem}{\em (\cite[Teor.2.3]{Valdes})} Suppose that $G$ satisfies the hypothesis of the above Theorem. Then, the  $G$-structure is integrable if and only if  the tensors of torsion and curvature of the well-adapted connection vanishes.
\end{theorem}

\subsection{Some important connections associated to an almost bi-Lagrangian manifold}

Now, we shall focus our attention on almost bi-Lagrangian manifolds.  We know (see Corollary \ref{existen}) the existence of the well-adapted connection.  Remembering Theorem \ref{equivalencia} we can obtain two equivalent versions of the following theorem, which allows us to have an explicit expression of the well-adapted connection.: 

\begin{theorem} {\em \cite[Cor. 3.1 and Theor. 3.1]{ES1999} }

 {\rm (bi-Lagrangian version).} An almost bi-Lagrangian manifold $(M,\omega,D_1, D_2)$ admits
functorial connections and its well adapted connection is the
unique linear connection $\nabla $ verifying $ \nabla \omega = 0, \quad \nabla F=0$ and
\[  \omega (F({\rm Tor}_{\nabla } (X , \pi _{1} Y)) , \pi _{2} Z)
 - \omega (F({\rm Tor}_{\nabla } (X , \pi _{2} Z)) , \pi _{1}Y) = 0 ,\]
for every vector fields $ X ,Y,Z$ on $M$, where $F^{+}=\frac{I+F}{2},\, \pi _{2}=\frac{I-F}{2}$ are the projections
over $D_1$ and $D_2$.
 
\bigskip

{\rm  (para-Hermitian version).} An almost para-Hermitian manifold $(M,g,F)$ admits functorial connections and its well adapted connection is the unique linear connection $\nabla $ verifying $ \nabla g = 0, \quad \nabla F=0$ and
\[  g ({\rm Tor}_{\nabla } (X , \pi _{1} Y) , \pi _{2} Z)
 - g ({\rm Tor}_{\nabla } (X , \pi _{2} Z) , \pi _{1}Y) = 0 ,\]
for every vector fields $ X ,Y,Z$ on $M$, where $\pi _{1}=\frac{I+F}{2},\, \pi _{2}=\frac{I-F}{2}$ are the projections over $F^+$ and $F^-$.

\end{theorem}

Using the symplectic version of the above theorem, observe that in the integrable case, i.e., $N_{F}=0$, then the 
well-adapted connection is  torsionless, thus proving that it coincides with the bi-Lagrangian connection.

On the other hand, one can define another linear connection on an almost para-Hermitian manifold (introduced by Libermann \cite{Lb}  in 1954 and characterized by Cruceanu and Etayo \cite{CE}):

\begin{theorem} {\em{\cite[Prop. 3.1]{CE}}} Let $(M, g , F)$ be an almost para-Hermitian manifold. There exits a unique linear connection   $\tilde{\nabla }$
satisfies the conditions
\newcounter{Item1000}\begin{list}{ \roman{Item1000})} {\usecounter{Item1000} }

\item  $\tilde{\nabla } F =0$.
\item $\tilde{\nabla } g =0$.
\item ${\rm Tor}_{\tilde{\nabla }} (X_1 , X_2)= 0, \quad \forall X_1 \in F^+ (M) , X_2 \in F^- (M)$.
\end{list}
This connection will be named the {\em Libermann connection} of  $(M,g,F)$.
\end{theorem}

One also can define the {\em natural connection}, see \cite[Proposition 5.2]{CE}, which is given by $\nabla _{X}= \frac{1}{2}(\nabla ^{g}_{X} + F\circ \nabla ^{g}_{X}  \circ F)$, where $\nabla ^{g}$ denotes the Levi-Civita connection of the neutral metric $g$. Observe that $F\circ \nabla ^{g}_{X}  \circ F$ is also a linear connection parallelizing $F$ and $g$. Similar connections have been defined on special complex manifolds (see \cite[Proposition 3]{ACD}).

Taking into account  Theorems \cite[Prop. 4.3]{ES1999} and \cite[Proposition 5.3]{CE} and Theorems \ref{metrica} and \ref{equivalencia} of the present paper, we obtain:

\begin{theorem}  (a) If $M$ is an almost bi-Lagrangian manifold then the almost bi-Lagrangian  and the Libermann connections coincide.

(b) If $M$ is a  bi-Lagrangian manifold then the well adapted, the Libermann, the natural, the  bi-Lagrangian, the almost
bi-Lagrangian  and the
Levi-Civita connections coincide.
\end{theorem}

Some characterizations of 1-para-Hermitian and 1-para-K\"{a}hler manifolds have been obtained in \cite{CE} by using selected connections. Moreover, the set of all the connections parallelizing  $F$ and $g$ of an almost Hermitian manifold is also obtained.

\section{The Holonomy of the leaves}

Let us consider a bi-Lagrangian manifold $M$. Then, the metric and the symplectic form vanish when are restricted to the leaves of both foliations.
Let $L$ be any leaf of one of  the foliations ${\cal F}_{i},\, i=1,2$. As $\nabla _{X}Y\in TL$ for every vector fields tangent to the leaf $L$ one obtains that $L$ is a totally geodesic submanifold of $M$. This make sense even though $L$ is
$g$-isotropic and one cannot decompose $T_{x}M$ as a direct sum of $T_{x}L$ and its orthogonal complement when $x\in L$. It means that the parallel transport, with respect to $\nabla$, along curves contained in the leaf $L$ moves vectors tangent to $L$ to vectors tangent to $L$, or, equivalently, the geodesics of $(M,\nabla )$ with initial point and derivative in $L$ are contained in $L$ (cfr. \cite[vol. II, pages 54-59]{KN}). Then, the connection $\nabla$ can be restricted to any leaf $L$, although $L$ is a Lagrangian $g$-isotropic manifold.

Now, we consider a point $p\in M$ and the leaves $L_{1}$ and
$L_{2}$, one of each foliation, through the point $p$. In this
section we shall obtain information about the number $N(p)$  of
points in the intersection $L_{1}\cap L_{2}$.  If $L_{1}\cap
L_{2}=\{ p\}$, i.e., if $N(p)=1$, we shall say that $p$ has the
{\em  trivial intersection property}. Observe that one can ask
about the number $N(p)$ when one has a manifold with two
transversal foliations.

\subsection{Bi-Lagrangian surfaces}
\label{surfaces}
The geometry of a bi-Lagrangian surface is quite simple, because every almost symplectic form is closed and every
1-dimensional distribution is involutive and Lagrangian. We shall show enough examples to prove that there are no relations among the following concepts: compactness, flatness and trivial intersection property.

If $M$ is a bi-Lagrangian surface, then $M$ is an orientable Lorentz surface. On the other hand, one has:

\begin{proposition}  Let $(M,g)$ be an orientable  Lorentz surface. Then $M$ is a bi-Lagrangian surface.
\label{lorentz}
\end{proposition}

{\em Proof.} We  follow the idea of \cite[Proposition 2.1]{ES2000}. As $(M,g)$ is a Lorentz surface then the nullcone at any point is given by two straight lines. Moreover, as $(M,g)$ is a orientable we can numerate both lines of the nullcone, thus definining two 1-dimensional distributions $V_{1}$ and $V_{2}$, satisfying the following properties, for all $x\in M$:
(1) $T_{x}M=(V_{1})_{x}\oplus (V_{2})_{x}$; 
(2) $dim\: (V_{1})_{x} = dim\: (V_{2})_{x}=1$; 
(3) $g\mid _{(V_{1})_{x}}=g\mid _{(V_{2})_{x}}=0$. 
Then, one can define an almost product structure $F$ such that the
$(V_{1})_{x}$ (resp. $(V_{2})_{x}$) is the eigenspace associated to +1 (resp. -1), for each $x\in M$. Then, by a straightforward computation, one concludes that $(M,F,g)$ is an almost para-Hermitian manifold, and  one defines the almost symplectic form $\omega$ given by $\omega (X,Y)=g(FX,Y)$ thus proving that $(M,F,g)$ is a parak\"{a}hlerian
(=bi-Lagrangian) surface $\Box$

\begin{remark}   One can define a Lorentzian metric over a non-orientable manifold (see, e.g., \cite[page 145, figure 5]{O'N}, where a time-orientable metric is defined on the M\"{o}bius band, in this case the lines of the nullcone do not define two distributions, as one can easily check). Moreover, one can define a Lorentzian metric on every non-compact manifold \cite[page 149]{O'N}, but in the compact case one only can define a Lorentz metric if the manifold has Euler number $\chi (M)=0$. Thus,  the unique closed (i.e., orientable and compact) surface which admits a Lorentz metric is the torus.
\label{superflorentz}
\end{remark}

Then we shall obtain the following examples:

\begin{center}
\begin{tabular}{|c|c|c|c|l|}
\hline
compact & flat & trivial & no. ex. & examples\\
\hline
yes & yes & yes &  1& {\em A flat torus.}\\
\hline
yes  &yes & no  & 2  &{\em Other flat tori.} \\
\hline
yes & no & yes &  3 &{\em The Clifton-Pohl torus.}\\
\hline
yes & no & no &  4 & {\em The Clifton-Pohl torus.}\\
\hline
no & yes & yes &  5&{\em The Minkowski plane.}\\
\hline
no  &yes & no  &   6 &  {\em The Lorentz cylinder.}\\
\hline
no & no & yes & 7a& {\em The ruled hyperboloid.} \\

   &        &         &     7b    & {\em The punctured plane.}\\

 &   &   &   7c& {\em The Schwarzschild half-plane.}\\
\hline
no & no & no & 8 & {\em A Clifton-Pohl cylinder.}\\
\hline
\end{tabular}
\end{center}
\vspace{5mm}

We shall also answer the following two questions:

\begin{itemize}

\item If $M$ is a bi-Lagrangian manifold, does every leaf of one foliation intersect every leaf of the other one? We shall show that the answer is negative, obtaining  (see Example 7a) a bi-Lagrangian surface with two leaves, one of each foliation, which do not intersect between them.

\item If $M$ is a bi-Lagrangian manifold, is the number of points of intersection of both leaves through a point $p$ independent of this point $p$? We shall show (see Examples 3 and 4) that the number $N(p)$ depends on the considered  point $p$.

\end{itemize}

\bigskip
{\bf Examples 5, 1 and 2:}

 Let us consider the plane $\mathbb{R}^{2}$, with global coordinates $(x,y)$ endowed with the symplectic form $\omega =\left( \begin{array}{cc} 0 & 1 \\ -1 & 0\end{array} \right) $. Then, any straight line of $\mathbb{R}^{2}$ defines a Lagrangian foliation.
Let us consider now the bi-Lagrangian structure determined by the diagonal lines $\{ x-y=0\}$ and $\{ x+y=0\}$. The almost product structure $F$ attached to these foliations is given by the matrix $F=\left( \begin{array}{cc} 0 & 1 \\ 1 & 0\end{array} \right) $ and the matrix expression of $g$ is $g=\left( \begin{array}{cc} -1 & 0 \\ 0 & 1\end{array} \right) $, which is a Lorentzian metric in $\mathbb{R}^{2}$. (Equivalently, we could consider the {\em Minkowski plane}
$\mathbb{R}^{2}_{1}=(\mathbb{R}^{2},g)$, whose nullcone is defined by the diagonal lines; taking into account Proposition \ref{lorentz}, we obtain a
bi-Lagrangian structure).

\bigskip

Let us consider the flat torus $\mathbb{T}^{2}=\mathbb{R}^{2}/\mathbb{Z}^{2}$ having the unit square of vertices $(\pm \frac{1}{2},\pm \frac{1}{2})$ as fundamental region. Then a leaf of one foliation meets a leaf of the other one in exactly two points. This is an example of compact flat bi-Lagrangian surface with $N(p)=2$, for all $p\in \mathbb{T}^{2}$.

Let us consider now another flat torus $\mathbb{T}^{2}=\mathbb{R}^{2}/\mathbb{Z}^{2}$ having the unit square of vertices
$(\pm \frac{1}{\sqrt{2}},0)$ and $(0,\pm \frac{1}{\sqrt{2}})$ as fundamental region, i.e., the  square in the above example rotated (in Riemannian sense) an angle of $\frac{\pi }{4}$. Then a leaf of one foliation meets a leaf of the other one in exactly one point.

In the above two examples any leaf of any foliation defines a
torus knot, because its slope with respect to the lattice defined
by the square is a rational number. If one rotates the square in
such a way that the Lagrangian lines $\{ x-y=0\}$ and $\{ x+y=0\}$
have rational slope one obtains torus knots, which have a finite
set of points of intersection. If the slope is not a rational
number, then the intersection is an infinite set.

\bigskip

{\bf Example 6}

 Let us consider the {\em Lorentz cylinder} $M$ as is defined in \cite[page 148]{O'N}. One can consider $M=S^{1}_{1}\times \mathbb{R}^{1}$ viewed as a cylinder in $\mathbb{R}^{3}_{2}$. Its geodesics coincide with ones of the Riemannian standard cylinder. Null geodesics can be parameterized as $(\pm \cos s,\sin s,s+c)$. Obviously both null geodesics through a point of $M$ intersect between them in an infinite set of points. Taking into account Proposition \ref{lorentz} we can conclude that the Lorentz cylinder is a non-compact flat bi-Lagrangian surface with non-trivial intersection.

\bigskip

{\bf Examples 7a, 7b and 7c}

Now we shall show three examples of  non-compact non-flat bi-Lagrangian surfaces such that leaves of both foliation meet in a unique point: the ruled hyperboloid, the punctured plane and the Schwarzschild half-plane. In particular, we shall show that a bi-Lagrangian manifold may admit two leaves, one of each foliation, which do not intersect between them (in the ruled hyperboloid one can take two parallel straight lines in ``antipodal" points of the equator).

\bigskip

{\bf a)} The following idea is due to Bejan in (see \cite[page 26]{Be}). Let us consider the Lorentzian space $\mathbb{R}^{3}_{1}$, an orthonormal basis $\{ e_{1},e_{2},e_{3}\}$ verifying
$-g(e_{1},e_{1})=g(e_{2},e_{2})=g(e_{3},e_{3})=1$, and the pseudosphere (in the sense of O'Neill's book \cite [page 110 ]{O'N}) $S^{2}_{1}(r)=\{ -x^{2}_{1}+x^{2}_{2}+x^{2}_{3}=r^{2}\}$. Then, $S^{2}_{1}(r)$ is a Lorentzian surface of positive curvature $\frac{1}{r^{2}}$ when one consider the induced metric $g|_{S^{2}_{1}(r)}$. As a quadric surface, $S^{2}_{1}(r)$ is a twofold {\em ruled hyperboloid}.

One can easily prove that $ \{ X_{1},X_{2}\}$ is an orthonormal basis of the tangent plane to $S^{2}_{1}(r)$ verifying
$-g( X_{1}, X_{1})=g( X_{2}, X_{2})=1$, when one defines:

$$ X_{1}= \frac{1}{rf}(f^{2}e_{1}+x_{1}x_{2}e_{2}+x_{1}x_{3}e_{3}) \hspace{3mm}; \hspace{3mm}
X_{2}= \frac{1}{f}(-x_{3}e_{2}+x_{2}e_{3})$$

\noindent with $f=\sqrt{r^{2}+x_{1}^{2}}$. Then, let $F$ be the (1,1)-type tensor field defined by $F(X_{1})=X_{2}$ and $F(X_{2})=X_{1}$, which has $Y=X_{1}+X_{2}$ and $Z=X_{1}-X_{2}$ as the eigenvectors fields associated to the eigenvalues +1 and -1.

Then, $(S^{2}_{1}(r),g,F)$ is a para-K\"{a}hlerian manifold with isotropic distributions generated by $Y$ and $Z$. A direct calculation shows that the integral curves of both fields $Y$ and $Z$ are the straight lines of the hyperboloid $S^{2}_{1}(r)$. In fact, this hyperbolid is the  bi-Lagrangian symmetric space $SL(2,\mathbb{R})/\mathbb{R}^{*}$ (see, e.g., \cite[p. 533]{K2}).

\vspace{5mm}

{\bf b)} Let us consider the {\em punctured plane} $\mathbb{R}^{2}-\{ (0,0)\}$ with global coordinates $(x,y)$ and the metric $g=\frac{2}{x^{2}+y^{2}}\, dx\otimes dy$. Then $(\mathbb{R}^{2}-\{ (0,0)\} ,g)$ is a Lorentzian surface, with non-constant curvature $K$, $-2\leq K\leq 2$, as one can easily check, using \cite[exercise 8, page 156]{O'N}. The nullcone at a point  is given by the horizontal and vertical lines through the point and then, this surface satisfies the hypothesis of Proposition \ref{lorentz}, thus providing an example of non-compact non-flat
bi-Lagrangian surface  such that the intersection of the leaves through a point reduces to the point.

Observe that if $p=(a,0)$ with $a>0$, then the horizontal leaf through is the positive $x$-axis, which does not intersect the vertical leaves $\{ x=b\}$ when $b<0$, thus proving that there exists a leaf that does not intersect an infinity of leaves of the other foliation.

 We shall use this example on the following examples 3 and 4: the homotheties of centre $(0,0)$ are isometries of the punctured plane $(\mathbb{R}^{2}-\{ (0,0)\} ,g)$; this property allows to define a Lorentzian metric in the torus in such a way that $\mathbb{R}^{2}-\{ (0,0)\}$ is an isometric covering. This manifold is named the {\em Clifton-Pohl torus}.

\vspace{5mm}

{\bf c)} The {\em Schwarzshild half-plane} is defined in, e.g., \cite[pag. 152-152]{O'N}). For a constant $M>0$ let $h(r)=1-(2M/r)$ and $P=\{ (t,r)\in \mathbb{R}^{2}, r>2M\}$, endowed with the Lorentzian metric $g=-h dt\otimes dt +h^{-1}dr\otimes dr$. Then $(P,g)$ is a surface of constant curvature $2M/r^{3}\, >0$. The null geodesics are obtained in \cite{O'N} (see, in particular, figure 7 in page 153, where it is shown that a geodesic of one family intersects in one point exactly any geodesic of the another family). Then, $(P,g)$ is a non-compact non-flat surface, and, taking into account Proposition \ref{lorentz}, $(P,g)$  is a bi-Lagrangian manifold.

\bigskip

{\bf Examples 3 and 4}

We shall show that in the Clifton-Pohl torus there exist two kinds of leaves: compact and non-compact. We shall show:

a)
Thus, as a consequence, we have obtained an example of bi-Lagrangian manifold on which the function $N:M\to \mathbb{N}\cup \{ \infty \}$ is not constant.
\bigskip

The definition of the {\em Clifton-Pohl torus} follows from the above example 7b (see, e.g., \cite[page 193]{O'N}): As homotheties are isometries of $(\mathbb{R}^{2}-\{ (0,0)\} ,g)$, one can consider the group $\Gamma =\{ \mu ^{n}\}$ generated by the homothety $\mu (x,y)=(2x,2y)$. $\Gamma$ is properly discontinuous, and  $T=M/\Gamma $ is a Lorentzian surface. Topologically $T$ is the closed annulus $\{ 1\leq \sqrt{x^{2}+y^{2}}\leq 2\}$, with boundaries points identified under $\mu$, i.e., it is a torus, named the Clifton-Pohl torus. The four compact leaves are the circles obtained intersecting the coordinates axes with the above annulus. Any other leaf is topologically a real line, which accumulates over two of the above circles, as one can easily see. Then properties a) and b) above are obvious. Finally, observe that there exist leaves, one of each foliation, without intersection: the compact leaves.

\bigskip

{\bf Example 8}

The punctured plane, defined in example 7b, is a topological cylinder which covers the Clifton-Pohl torus  by a locally isometric submersion. From the point of view of the function $N$, both surfaces are quite different, because the punctured plane has trivial intersection whereas in the Clifton-Pohl torus there exist points with non-trivial intersection (indeed, almost every point has this property).

We shall define another topological cylinder $M$, which also covers the Clifton-Pohl torus by a locally isometric submersion, but in such a way that it preserves the non-triviality property. The idea looks like the construction of the Riemann surface associated to the complex logarithm. Let us consider as fundamental region the annulus $\{ 1\leq \sqrt{x^{2}+y^{2}}\leq 2\}$, with boundary points identified under $\mu$, as in the above example (5). Let us consider left and right semi-annulus, which are obtained cutting the above one by the $y$-axis:

$$L=\{ 1\leq \sqrt{x^{2}+y^{2}}\leq 2, x\leq 0\}\hspace{3mm};\hspace{3mm}R=\{ 1\leq \sqrt{x^{2}+y^{2}}\leq 2, x\geq 0\}$$

We define a countable family $\{ L_{n}\} _{n\in \mathbb{Z}}$, where $L_{n}$ is $L$, and another family $\{ R_{n}\} _{n\in \mathbb{Z}}$, where $R_{n}$ is $R$. Finally we identify the subset $\{ x=0, y>0\}$ of $L_{n}$ with the subset $\{ x=0, y>0\}$ of $R_{n}$ and the subset $\{ x=0, y<0\}$ of $R_{n}$ with the subset $\{ x=0, y<0\}$ of $L_{n+1}$. Then, the family $M=(\bigcup L_{n})\cup (\bigcup R_{n})$
with the identification topology (and identifying the boundaries points as in the Clifton-Pohl torus) defines a topological cylinder, which has a bi-Lagrangian structure obtained by lifting that of the Clifton-Pohl torus, and satisfies the desired conditions, as one can easily show.

\subsection{Higher dimensional bi-Lagrangian manifolds}

The same problem is  completely different for higher dimensional manifolds, becasuse neutral $\neq$ Lorentz if ${\rm dim}(M)>2$. Of course, one can easily define in $\mathbb{R}^{2n}$ a bi-Lagrangian structure with $N(p)=1$, for all point $p$, thus obtaining a non-compact flat
bi-Lagrangian manifold with trivial intersection. 

In order to obtain some interesting examples we can use tangent and cotangent  bundles (in the non-compact case) and the square of a Lie group $G\times G$ and the product of two bi-Lagrangian manifolds. We begin with the table of examples:

\begin{center}
\begin{tabular}{|c|c|c|c|l|}
\hline
compact & flat & trivial & no. ex. & examples\\
\hline
yes & yes & yes &  1& {\em A flat torus $\mathbb{T}^{2n}$.}\\
\hline
yes  &yes & no  & 2  & $\mathbb{T}^{2}_{1}\times \mathbb{T}^{2}_{2}$,  $\mathbb{T}^{2}_{1}$ {\em (resp. $\mathbb{T}^{2}_{2}$ ) with $N(p)=1$ (resp.  $N(p)=2$)} \\
\hline
yes & no & yes &  3 &{\em $G\times G$, $G$ being a non-flat torus $\mathbb{T}^{n}$.}\\
\hline
yes & no & no &  4 & {\em $T\times T$, $T$ being a Clifton-Pohl torus}  \\
\hline
no & yes & yes &  5a&{\em The neutral Euclidean space.}\\
   &        &         &     5b    &{\em $G\times G$, $G$ being a flat cylinder $\mathbb{T}^{n}\times \mathbb{R}^{m}$.}\\
   &        &         &     5c    & {\em $T(\mathbb{T}^{n})$.}\\
\hline
no  &yes & no  &   6 &  {\em $T(M)$, $M$ the Hantzsche-Went manifold}\\
\hline
 no & no & yes  &     7a    & {\em A non-flat cylinder $\mathbb{T}^{4}\times \mathbb{R}^{2k}$. }\\
   &        &         & 7b& {\em Kaneyuki examples.} \\

   &        &         &     7c    & {\em The paracomplex projective space.}\\
\hline
no & no & no & 8 & $\mathbb{R}^{2}\times T$, {\em $T$ being a Clifton-Pohl torus}\\
\hline
\end{tabular}

\end{center}
\vspace{5mm}

\subsubsection{The square of a  Lie group} 

We show that one can obtain examples in the compact case taking a compact Lie group $G\times G$, where $G$ is a compact Lie group endowed with a left invariant Riemannian metric $<\, ,\, >$. Moreover, one can also use this construction in the non-compact case.

Let $G$ be a Lie group and let $<\, ,\, >$ be an inner product in the tangent space at the identity, $T_{e}G$. Then one can define a left invariant Riemannian metric on $G$ by means of 

$$<v_{a} ,w_{a} >=< (L_{a^{-1}})_{*}(v_{a}) ,(L_{a^{-1}})_{*}(w_{a})  >$$

\noindent for all $v_{a},w_{a}\in T_{a}G$.

Let us consider the product Lie group $G\times G$, which has two transversal foliations: if $(a,b)\in G\times G$, the leaves  through $(a,b)$ are $\{ a\} \times G$ and $G\times \{ b\} $. Taking into account the isomorphism $T_{a}G\times T_{b}G\thickapprox T_{(a,b)}G\times G$ we can define: 

\bigskip

$$ g\left( (v_{a},v_{b}),(w_{a},w_{b}) \right)  \\   \\
  =\,  < (L_{a^{-1}})_{*}(v_{a}) ,(L_{b^{-1}})_{*}(w_{b})>
+ <(L_{a^{-1}})_{*}(w_{a}) , (L_{b^{-1}})_{*}(v_{b}) >. $$

\bigskip

One easily check the following properties: 

1)  $  g\left( (v_{a},0),(w_{a},0) \right) =0=g\left( (0,v_{b}),(0,w_{b}) \right) $, which shows that the leaves are $g$-isotropic. 

2)  $g$ is a neutral metic on $G\times G$: if $\{  e_{1}, \ldots ,e_{n}\} $ is a $<\, ,\, >$-orthonormal  basis of $T_{e}G$, then
$\{ (L_{a*}(e_{1}),0), \ldots , (L_{a*}(e_{n}),0), (0,L_{b*}(e_{1})), \ldots ,(0,L_{b*}(e_{n})) \} $ is a basis of $T_{(a,b)}G\times G$ such that $g$ has matrix

$$ \left(  \begin{array}{cc} 0 & I \\ I & 0 \end{array} \right)  $$

\noindent respect to it, thus proving that $g$ is a neutral metric of signature $(n,n)$. If  $\Delta :G\to G\times G$ denotes  the diagonal embedding then one has $g|_{\Delta (G)}=2<\, ,\, >$, which proves that $(G\times G,g)$ is non-flat if $(G,<\, , \, >)$ is non-flat.

3) Let $F$ be the almost product structure associated to the foliations ${\cal F} _{1}=\left\{ \{ a\} \times G, a\in G\right\}$  and ${\cal F} _{2}=\left\{ \{ G\times \{ b\},b\in G\right\} $. Then, the almost symplectic form $\omega $ given by $\omega (X,Y)=g(FX,Y)$ vanishes on the leaves of both foliations. Then,  $(G\times G,\omega ,{\cal F} _{1},{\cal F} _{2})$ is a
 para-Hermitian manifold.

\bigskip 

In order to prove that $(G\times G,\omega ,{\cal F} _{1},{\cal F} _{2})$ is a bi-Lagrangian manifold we have to show that $\omega$ is closed. As $\omega $ is a 2-form we have:

\begin{eqnarray*}
d\omega (X,Y,Z)& =&  X(\omega (Y,Z)) + Y(\omega (Z,X)) + Z(\omega (X,Y)) \\
 & - & \omega ([X,Y],Z)- \omega ([Y,Z],X)-\omega ([Z,X],Y)   
\end{eqnarray*}

As is well known, the Lie algebra of $G\times G$ is the Lie algebra product $\mathfrak{g}\times \mathfrak{g}$ given by

$$ [(X,Y),(X',Y')]=([X,X'],[Y,Y'])             $$

Observe that  if $(X,0),(X',0)\in T({\cal F} _{2})$, then  $ [(X,0),(X',0)]=([X,X'],0)$, which is a vector field of the same distribution $T({\cal F} _{2})$. The same is true for  ${\cal F} _{1}$. 

Let us consider the global basis of vector fields given in 2) above:

$$\{ (L_{a*}(e_{1}),0), \ldots , (L_{a*}(e_{n}),0), (0,L_{b*}(e_{1})), \ldots ,(0,L_{b*}(e_{n}))\} $$

\bigskip

Then $\omega ((L_{a*}(e_{i}),0)  ,(L_{a*}(e_{j}),0)  )=0= \omega (  (0,L_{b*}(e_{i}) , (0,L_{b*}(e_{j}) )$, because the foliations are Lagrangian and 

\begin{eqnarray*}
X(\omega ((L_{a*}(e_{i},0)  , (0,L_{b*}(e_{j})  )   & =& X(g(-(L_{a*}(e_{i},0), (0,L_{b*}(e_{j}))\\
&  = & X(<-e_{i},e_{j}>)=X(\delta _{ij})=0  
\end{eqnarray*}

Thus, we have proved that the three first terms of $d\omega (X,Y,Z)$ vanish when $X,Y,Z$ are vector fields of our basis. For the other terms, we have to study the Lie brackets. As the distributions are involutive and Lagrangian, all the three last terms vanish except to those of the form:

\begin{eqnarray*}
\omega ([(L_{a*}(e_{i},0) , (L_{a*}(e_{j},0) ]   , (0,L_{b*}(e_{k})   )= \omega (([L_{a*}(e_{i}),L_{a*}(e_{j})],0),(0,L_{b*}(e_{k})   )=\\
\omega ((L_{a*}[e_{i},e_{j}],0),(0,L_{b*}(e_{k})   )=   <-[e_{i},e_{j}],e_{k}>=<-c_{ij}^{\ell }e_{\ell },e_{k}>=-c_{ij}^{k}
\end{eqnarray*}

\noindent where $c_{ij}^{\ell }$ are the structure constants associated to the basis $\{ e_{1}, \ldots ,e_{n}\} $ of $T_{e}G$. In the case when $G$ is an abelian group, the structure constants vanish. Thus, we have proved the following result:

\begin{proposition}
\label{grupos}

Let $G$ be an abelian Lie group. Then $G\times G$ can be endowed with a bi-Lagrangian structure. Moreover, $N(a,b)=1$, for all $(a,b)\in G\times G$. 

\end{proposition}

And then we have:

\bigskip
{\bf Examples 1 and 3}

Let us consider $G$ an abelian connected compact Lie group. As is well known, then $G=\mathbb{T}^{n}$ is a 
n-dimensional  torus. Then $G\times G$, which is also a torus $\mathbb{T}^{2n}$, is a compact bi-Lagrangian manifold  with the trivial intersection property. If $<\, , \, >$ is a non-flat metric on $G$ then $G\times G$ is a non-flat manifold. For example, if  $G=\mathbb{T}^{2}\subset \mathbb{R}^{3}$ is the standard non-flat torus.

\subsubsection{Product of bi-Lagrangian manifolds} 

We shall show that the product of two bi-Lagrangian manifolds is also a bi-Lagrangian manifold. We shall use the 
para-K\"{a}hler terminology.

\begin{proposition}
\label{producto}
The product of two bi-Lagrangian manifolds is also a bi-Lagrangian manifold.
\end{proposition} 

{\em Sketch of Proof.}  Let us denote by $(M,F^{+},F^{-},g_{M}))$ and $(N,P^{+},P^{-},g_{N}))$ the both para-K\"{a}hler manifolds. Then $(M\times N,g)$ with $g=\pi ^{*}(g_{M})+\sigma ^{*}(g_{N})$ is a neutral manifold, $\pi :M\times N\to M$ and $\sigma :M\times N\to N$ being the projections (cfr., e.g.,\cite[p. 57]{O'N}). One can define an almost product structure $H$ on $M$ satisfying $H^{+}=F^{+}\oplus P^{+}$ and $H^{-}=F^{-}\oplus P^{-}$. Taking into account that the Levi-Civita connection of $g_{M}$ (resp. $g_{N}$) parallelizes $F$ (resp. P) and the following relations:

$\nabla _{X}Y= \nabla ^{M}_{X}Y$; for all $X,Y$ tangent to $M$;

$\nabla _{V}W= \nabla ^{N}_{V}W;$ for all $V,W$ tangent to $N$;

$\nabla _{X}W= \nabla _{V}Y=0$; for all $X,Y$ tangent to $M$ and $V,W$ tangent to $N$; 

\noindent one easily checks that the Levi-Civita connection of $g$ parallelizes $H$ thus finishing the proof $\Box$

\begin{remark}
\label{observproducto}

In the above proposition both manifolds $M$ and $N$ can be of different dimension. If any of both manifolds is non-flat, then $M\times N$ is non flat. If both manifolds are flat, $M\times N$ is also flat. In any case, the sectional curvature of a plane spanned by two vectors, each one tangent to each manifold, vanishes (cfr., e.g., {\rm \cite[p. 89]{O'N}}).

\end{remark}

\bigskip
{\bf Example 2}

Let us consider the examples 1 and 2 in Section \ref{surfaces}:  we define the product $\mathbb{T}^{2}_{1}\times \mathbb{T}^{2}_{2}$ of two torus,  $\mathbb{T}^{2}_{1}$ (resp. $\mathbb{T}^{2}_{2}$ ) with $N(p)=1$ (resp.  $N(p)=2$). Then the product manifold has no trivial intersection property. 

\bigskip
{\bf Example 4}

The same idea runs for the Clifton-Pohl torus $T$: the product $T\times T$ is a compact non-flat manifold with no trivial intersection property.

\bigskip
{\bf Example 8}

If we change a torus by a plane we have $\mathbb{R}^{2}\times T$ ($T$ being a Clifton-Pohl torus), which is no compact 
non-flat with no trivial intersection property.

\subsubsection{Holonomy of a flat Riemannian manifold }
\label{HW}

We need to study some gometric properties of flat Riemannian manifold and tangent bundles in order to obtain significative examples in the non-compact case.

Let $(M,\cal{F})$ be a manifold endowed with a foliation. The {\em holonomy} of $\cal{F}$ measures the intersection of any leaf of $\cal{F}$ with a transverse submanifold. In a general situation this transversal submanifold can intersect the leaf in a complicated set, but in the case of a flat bundle the intersection of the fibres with the horizontal leaves is a discrete set. In particular,  the tangent bundle of a flat Riemannian manifold has discrete holonomy and this topological notion of holonomy of the horizontal foliation coincides with that of geometric holonomy obtained by parallel transport (see, e.g. \cite[chapter 4]{T}). Moreover, if $(M,g)$ is a flat Riemannian manifold, then there exists a canonical map $\pi _{1}(M)\to {\rm Hol}(M)$, from the fundamental group of $M$  onto the holonomy group . If $M$ is simply connected, then its holonomy is trivial. The non-simply connected case  explains the  Aharonov-Bohm effect.

The following results are obtained and quoted in \cite{RT}. If $(M,g)$ is a compact flat Riemannian manifold of dimension $n$, then $M$ is a quotient manifold of $\mathbb{R}^{n}$ by a torsion-free discrete subgroup $\Gamma$ of the group of isometries $I(\mathbb{R}^{n})=O(n)\ltimes \mathbb{R}^{n}$, $\Gamma $ is isomorphic to $\pi _{1}(M)$ (because $\mathbb{R}^{n}$ is the universal covering of $M$), the holonomy group ${\rm Hol}(M)$ is a finite subgroup of $O(n)$  and one has an exact sequence

$$ 0 \to \Lambda \to \pi _{1}(M) \to {\rm Hol}(M)\to 1 $$

Now we are interested on  compact flat Riemannian manifolds with non-trivial holonomy. Hantzsche and Went obtained in 1935 the unique example of a 3-dimensional compact flat manifold with first Betti number zero. The group of holonomy  of the {\em Hantzsche-Went manifold} is $\mathbb{Z}_{2}\oplus \mathbb{Z}_{2}$. In 1975 Cobb \cite{Co} obtained  an infinite   family of compact flat Riemannian manifolds of dimension $\geq 3$  with first Betti number zero and  with holonomy $\mathbb{Z}_{2}\oplus \mathbb{Z}_{2}$. This family has been considerably increased in \cite{RT}. In any case, for all $n\geq 3$ there exits at least a compact flat Riemannian manifold (with first Betti number zero) and with group of holonomy equal to $\mathbb{Z}_{2}\oplus \mathbb{Z}_{2}$.

On the other hand, compact flat Riemannian manifolds with holonomy $\mathbb{Z}_{2}^{k}$ have been obtained in \cite{DM}.

\subsubsection{The tangent bundle of a flat Riemannian manifold}
\label{tangentbundle}

We shall show: (1) the tangent bundle $TM$  of a flat Riemannian manifold $(M,g)$ admits a canonical
bi-Lagrangian structure; (2) its Lagrangian foliations are the vertical and the horizontal ones; (3) the canonical connection is flat, and then, the Lagrangian foliations are totally geodesic isotropic submanifolds of a flat bi-Lagrangian manifold.

Let $M$ be a $n$-dimensional manifold endowed with a metric $g$ and let
$\nabla$ be the Levi-Civita connection of $g$.
One can also introduce the  almost para-complex structure on $TM$ defined by
$ FX^V = - X^V , \: FX^H = X^H $, where $X$ is a vector field on $M$ and $V$ (resp. $H$) denotes the
vertical lift (resp. horizontal) to the tangent bundle (see \cite{YI}). (The opposite of this structure has been
introduced in \cite{C1}). Let $g^H$ be the horizontal lift of $g$ and $\nabla^H$ the horizontal lift of $\nabla$. Then,
 $(TM , F , g^H)$ is an almost para-Hermitian manifold and $\nabla^H g^H = 0$.

Let $\omega $ be the almost symplectic structure of $(TM , F , g^H)$, i.e., $\omega (\overline{X},\overline{Y})=g^{H}(F\overline{X},\overline{Y})$, for all vector fields $\overline{X},\overline{Y}\in \mathfrak{X}(TM)$.  One  obtains:

\[ \omega (X^V , Y^V ) = 0 =\omega (X^H , Y^H ) , \, \omega (X^V , Y^H ) = -
(g (X ,Y))^V =-\omega (X^H , Y^V ), \]

\noindent which proves that the vertical and horizontal distributions are Lagrangian.

 The connection $\nabla^H$ satisfies $ \nabla^H F = 0 , \quad \nabla^H g^H = 0$ (see \cite [Theorem 1]{E2}).
Taking account that $\nabla$ is symmetric,  one has $ {\rm Tor}_{\nabla^H} (X^V , Y^H) = ({\rm Tor}_{\nabla} (X,Y))^V = 0$,
for every $ X , Y$ vector fields on $M$, where ${\rm Tor}_{\nabla^H}$ is
the torsion tensor field of $\nabla^H$.

If $\nabla$ is symmetric  flat, i.e., ${\rm Tor}_{\nabla}=0$ and $R_{\nabla}=0$,
where $R_{\nabla}$ denotes the curvature tensor of $\nabla$, then $\nabla^H$
is also symmetric  flat, taking account the following result \cite[Propositions 7.3 and 7.4]{YI}:
Let  $\nabla$ be a symmetric connection. Then the connection $\nabla^H$ is
symmetric if and only if $R_{\nabla} =0$. In this case, one also has
$R_{\nabla^H}=0$.

Finally, taking into account all the above results in this section, we can conclude: if $(M,g)$ is  flat Riemannian manifold and  $\nabla$ is the Levi-Civita connection of $g$, then, the horizontal distribution is involutive, $(TM,\omega )$ is a
bi-Lagrangian manifold, $\nabla^H$ is the Levi-Civita connection of  $g^H$ and $\nabla^H$ is also the canonical connection of $(TM,\omega )$.

\subsubsection{The non-compact examples}

We shall study the examples 5, 6 and 7.

\bigskip

{\bf Examples 5a, 5b and 5c}

The neutral Euclidean space is the obvious example. Let $G$ be a flat cylinder $\mathbb{T}^{n}\times \mathbb{R}^{m}$. Then, by Proposition \ref{grupos} one has a flat cylinder $G\times G= \mathbb{T}^{2n}\times \mathbb{R}^{2m}$ which is a non-compact bi-Lagrangian manifold with the trivial intersection property.

Let $(M,g)$ be a simply connected flat Riemannian manifold and let $(TM,\omega )$ be the tangent bundle with the
bi-Lagrangian flat structure. The holonomy of $(M,g)$ is trivial and then vertical and horizontal leaves of $(TM,\omega )$ intersect in one point. This is the case of $M=\mathbb{R}^{n}$ endowed with the canonical Riemannian metric.

 The same property is true if $(M,g)$ is any  flat Riemannian manifold with trivial holonomy, even though it is not simply connected. For example, if $M$ is any flat torus $\mathbb{R}^{n}/\Gamma$ (in \cite{GHL} one can learn the basic properties of flat tori).

\bigskip

{\bf Example 6 }

We  show that there exists a family of manifolds of dimension $2n$, $\forall n\geq 3$, of non-compact flat manifolds with $N(p)=4$, for all point $p$.

 Let $(M,g)$ be a compact flat Riemannian manifold with holonomy $\mathbb{Z}_{2}\oplus \mathbb{Z}_{2}$ (see section \ref{HW}) and let us consider the tangent bundle endowed with the bi-Lagrangian structure $(TM,\omega )$ obtained in the above section \ref{tangentbundle}. Then vertical leaves (i.e., the fibres) and horizontal leaves intersect between them  in $4 =\sharp \{\mathbb{Z}_{2}\oplus \mathbb{Z}_{2}\}$ points. Observe that horizontal leaves are compact, whereas vertical ones are no compact, but all of them are totally geodesic submanifolds of a flat manifold. The same idea runs for compact flat Riemannian manifolds with another group of holonomy.

\bigskip

{\bf Examples 7a, 7b and 7c}

We shall begin with the non-flat cylinder. Let us consider the standard torus $\mathbb{T}^{2}\subset \mathbb{R}^{3}$, which is a non-flat Riemannian manifold, and let us consider the cylinder $G=\mathbb{T}^{2}\times \mathbb{R}^{k}$. Then, by using Proposition \ref{grupos}, the cylinder $G\times G$  can be endowed with a non-flat bi-Lagrangian structure having the trivial intersection property.

\bigskip

For the other examples, let us remember  a recent paper of Kaneyuki  where has studied bi-Lagrangian symmetric spaces, proving the following result:

\begin{proposition} {\em \cite[Lemma 2.1]{K2}} Let $(M=G/G_{0}, \omega ,{\cal F} _{1},{\cal F}_{2})$ be a bi-La\-grang\-ian symmetric space associated to a simple graded Lie algebra $\mathfrak{g}_{-1}+\mathfrak{g}_{0}+\mathfrak{g}_{1}$. Then, for any point $p\in M$ we have $N(p)=1$.

\end{proposition}

In this case, the proof is not difficult because one can use the exponential map of Lie groups. Taking into account the classification of bi-Lagrangian symmetric spaces of the above kind, obtained by Kaneyuki and Kozai \cite{KK}, one can find a lot of examples with the trivial intersection property. All of them are no compact, because such a manifold is always diffeomorphic to the cotangent bundle of a covering manifold of a certain Riemannian space.

  For example, the paracomplex projective space  
$$P_{n}(B) = Sl(n+1,\mathbb{R})/S(Gl_{0}(n,\mathbb{R})\times Gl_{0}(n,\mathbb{R}))$$
\noindent and the paraquaternionic projective space, cfr. \cite{E93},

$$P_{n,n}(\mathbb{C}) = Gl(n+1,\mathbb{C})/Gl(1,\mathbb{C})\times Gl(n,\mathbb{C})$$ 

\noindent are bi-Lagrangian symmetric spaces. We shall give an introduction to the geometry of the first space, which has been studied by several authors. In fact, this manifold is known as the {\em paracomplex projective space} because it is related with the paracomplex numbers (see, e.g., \cite{CFG}), and we shall consider its properties studied in \cite{GM1}, \cite{GM2}, \cite{GM3} and \cite{EF}. We shall denote it as $P_{n}(B)$.

Let $P_{n}(B)=\{ (u,v)\in \mathbb{R}^{n+1}\times  \mathbb{R}^{n+1}: <u,u>=<v,v>; <u,v>=1\}$ where $<\, ,\, >$ denotes the canonical metric. Then, $P_{n}(B)$ is a $2n$-dimensional manifold which is globally diffeomorphic to $TS^{n}$ by means of the map $(u,v)\mapsto (\frac{u+v}{\| u+v\| },u-v)$. This manifold admits a canonical almost product structure $F$ and a neutral metric $g$, making it a para-K\"{a}hlerian manifold, and then a bi-Lagrangian manifold, defining $\omega (X,Y)=g(FX,Y)$, for all vector fields $X,Y$ tangent to $P_{n}(B)$. In order to define these structures we need to introduce local coordinates.

Local charts $(U^{+}_{\alpha }, \psi _{\alpha })$ and $(U^{-}_{\alpha }, \psi _{\alpha })$ are defined on $P_{n}(B)$ by

$$U^{+}_{\alpha }=\{ (u,v)\in P_{n}(B): u^{\alpha }>0, v^{\alpha }>0\} ,\hspace{5mm} 
U^{-}_{\alpha }=\{ (u,v)\in P_{n}(B): u^{\alpha }<0, v^{\alpha }<0\} , $$

$$\psi _{\alpha }(u,v) = (\frac{u^{0}}{u^{\alpha }}, \ldots \widehat{\frac{u^{\alpha }}{u^{\alpha }}}\ldots \frac{u^{n}}{u^{\alpha }};
\frac{v^{0}}{v^{\alpha }}, \ldots \widehat{\frac{v^{\alpha }}{v^{\alpha }}}\ldots \frac{v^{n}}{v^{\alpha }})$$

\noindent where the hat \hspace{2mm} $\widehat{}$ \hspace{2mm} denotes a deleted element. The local coordinates are $(x^{i}=\frac{u^{i}}{u^{\alpha }}, y^{i}=\frac{v^{i}}{v^{\alpha }})$. In these local coordinates the para-K\"{a}hlerian structure is given by:

$$ F=\frac{\partial  }{\partial x^{i}}\otimes dx^{i} -\frac{\partial  }{\partial y^{i}}\otimes dy^{i}$$

$$ g= \sum_{i,j} \frac{2}{c(1+<x,y>)}\;  [dx^{i}\otimes dy^{i} -\frac{1}{1+<x,y>}x^{i}y^{j} (dy^{i}\otimes dx^{j}+dx^{j}\otimes dy^{i})]                         $$

\bigskip

Then,  $(P_{n}(B),F,g)$ is a para-K\"{a}hlerian space form of constant paraholomorphic sectional curvature $c$ (i.e., the planes $F$-invariant which are not $g$-degenerates have sectional curvature equal to $c$). The sectional curvature runs over all the real line $\mathbb{R}$, when one moves the planes (see also \cite{E3}). In local coordinates the Lagrangian foliations of the symplectic form $\omega$ defined as $\omega (X,Y)=g(FX,Y)$ are the eigenspaces $F^{+}$ and $F^{-}$ associated to the eigenvalues +1 and -1 of $F$, i.e., $F^{+}= \{  y^{1}=const.,\ldots ,y^{n}=const.\}$ and $F^{-}= \{  x^{1}=const., \ldots ,x^{n}=const.\}$, which meet in one point, in each chart.

Following \cite{GM3}, we consider the manifold $P_{n}(B)/\mathbb{Z}_{2}$, named the {\em reduced paracomplex projective space}, which is globally diffeomorphic to the tangent bundle of the real projective space, $P_{n}(B)/\mathbb{Z}_{2}\thickapprox TP_{n}(\mathbb{R})$. Then, $P_{n}(B)/\mathbb{Z}_{2}$ is also a para-K\"{a}hlerian space form of constant paraholomorphic sectional curvature $c$, because $P_{n}(B)$ is a paraholomorphically isometric twofold covering of $P_{n}(B)/\mathbb{Z}_{2}$.

On the other hand, in the case $n=1$ one has \cite{GM2} that $P_{1}(B)$  is a paraholomorphically isometric twofold covering of the ruled hyperboloid $H=\{ (x,y,z)\in \mathbb{R}^{3}: x^{2}+y^{2}-z^{2}=1/|c|\}$, whose para-K\"{a}hlerian structure is defined by the almost product structure determined by the straight lines and the metric induced by the
semi-Riemannian metric $\overline{g}=(c/|c|) (dx^{2}+dy^{2}-dz^{2})$ of $\mathbb{R}^{3}$. Taking into account both results, one has that $P_{1}(B)/\mathbb{Z}_{2}$ is paraholomorphically isometric to the above hyperboloid $H$.

\begin{remark}

 Following the ideas of the examples 7.a in section \ref{surfaces} and the present one,  a direct relation between pseudospheres and paracomplex projective spaces has been obtained in \cite{E3}: the
pseudosphere $S^{2n+1}_{n+1}$ is diffeomorphic to the product
$P_{n}(B)\times \mathbb{R}^{+}$. Moreover, one can define a
principal bundle $S^{2n+1}_{n+1}\to P_{n}(B)$ which allows to
obtain a Fubini-Study type metric on the paracomplex projective
space.
\end{remark}

\section{Open problems}

We want to end this work with a list of unsolved problems and open questions.

(1)There is no topological classification  for closed manifolds of dimension greater than three admitting local-product structures. For the three-dimensional case, such a manifold is homeomorphic to a Seifert manifold with zero Euler number (see \cite{Ma}). A similar open problem consists on obtaining a topological classification of closed manifolds admitting a
bi-Lagrangian structure. As we have seen in Remark \ref{superflorentz}, the torus is the unique closed surface admitting a bi-Lagrangian structure, but the problem remains open for higher dimensions.

If one considers para-K\"{a}hler space forms $(M,g,F)$, i.e.,
para-K\"{a}hler manifolds having constant metric sectional
curvature on the $F$-invariant planes, some results have been
obtained \cite{GM1995}: if $M$ is a complete and connected
manifold of dimension $2n>2$ and $c\neq 0$ then $M$ is
para-holomorphically isometric to a space $T(S^{n}/\Gamma )$,
where $\Gamma $ is a finite group with additional conditions. If
$n$ is even, then $M$ is paraholomorphically isometric to $TS^{n}$
or $TP_{n}(\mathbb{R})$ and $M$ is homogeneous. For the case
$c\neq 0$, ${\rm dim}M=2$ and for the cases $c=0$, ${\rm dim}M\geq
2$ the question remains open: a para-holomorphic classification
has not been found.

(2) Determine the group of paracomplex automorphisms and that of paracomplex isometries of a bi-Lagrangian manifold. In this topic the paper of Kaneyuki \cite{K2} must be the start point.

(3)The theory of real submanifolds of complex manifolds shows a large collection of interesting submanifolds, such as complex, totally real, Cauchy-Riemann, slant, generic, etc. In a recent paper \cite{E2003}, one of the authors have obtained results about the holomorphicness of a real submanifold of an almost Hermitian manifold. It would be interesting to obtain similar results about submanifolds of a symplectic manifold, or, at least, about submanifolds of a bi-Lagrangian manifold.

\vspace{5mm}

\noindent {\Large {\bf Acknowledgments}}

The work has been partially supported through grant DGICYT (Spain) (Project BFM2002-00141). The authors are also grateful to Professors V. Cruceanu, P. M. Gadea, J. Mu\~{n}oz Masqu\'{e}, R. F.  Picken, G. Thompson and S. Vacaru for their comments about some aspects of the topic.

\begin{footnotesize}

\end{footnotesize}

\end{document}